\documentclass[11pt]{article}
\usepackage{amsmath,amssymb}
\usepackage{graphicx,color}
\usepackage{bm}
\newcommand{\qed}{\hfill Q.E.D.\\}
\newtheorem{theorem}{Theorem}[section]
\newtheorem{proposition}[theorem]{Proposition}
\newtheorem{lemma}[theorem]{Lemma}

\newtheorem{remark}{Remark}

\newtheorem{definition}[theorem]{Definition}

\numberwithin{equation}{section}
\numberwithin{remark}{section}
\newcommand{\Proof}{{\bf Proof:}\quad}
\begin{document}
\title{Spectral problems about many-body Dirac operators\\ mentioned by Derezi\'{n}ski}
\author{Takashi Okaji\thanks{Department of Mathematics,
Graduate School of Science,Kyoto University}, 
Hubert Kalf\thanks{Mathematisches Institut der Universit\"{a}t M\"{u}nchen} and 
Osanobu Yamada\thanks{Department of Mathematical Sciences, Ritsumeikan University}}
\maketitle
\begin{abstract}
We consider spectral problems for many-body Dirac operators mentioned
by Derezi\'{n}ski in the IAMP News Bulletin of January 2012. In particular, 
we derive a representation of the Dirac Coulomb operator for a helium-like 
ion as a matrix operator of order sixteen. We show that it is essentially 
self-adjoint (under natural restrictions on the coupling constants), 
that the essential spectrum of its closure is the whole real line and that 
it has no eigenvalues.
\end{abstract}
\section{Introduction}
In \cite{D}, Derezi\'{n}ski mentioned open problems about many-body Dirac operators.
These problems were originally formulated by B. Jeziorski 
who is a chemist from the University of Warsaw (cf. J. Sucher \cite{S} p.6).
Among them there are spectral problems on Dirac-Coulomb operator $H_{DC}$ 
for a helium-like ion, which has the form
\begin{equation}
	H_{DC}=H(1,Z)+H(2,Z)+\frac{1}{| {\bf r}_{1}-{\bf r}_{2}|},
\end{equation} 
where 
\begin{equation}
	H(i,Z)=c\vec{\alpha}\vec{p}_{i}+mc^{2}\beta -\frac{Z}{|{\bf r}_{i}|} ,\ i=1,2
\end{equation} 
is the usual Dirac operator for an electron $i$ in the hydrogen-like ion
of charge $Z$ and of mass $m$.
In the above notation, ${\bf r}_{i}$ and $\vec{p}_{i}$, $i=1,2$  are 
a position vector and the momentum operator, respectively, of the $i$-th electron,
\[ 
	{\bf r}=(x_{1},x_{2},x_{3}),\ \vec{p}=-i\hbar\; \mathrm{grad}.
\] 
The vector $\vec{\alpha }$ is a vector operator whose components $\alpha _{1},\alpha _{2},\alpha _{3}$,
together with the operator $\beta \equiv \alpha _{4}$
are Hermitian matrices of order four satisfying the anti-commuting relations
\[ 
	\alpha _{j}\alpha _{k}+\alpha _{k}\alpha _{j}=2\delta _{jk}\
	(j,k=1,2,3,4).
\] 
Since the domain of the Dirac operator $H(i,Z)$ is a subspace of 
four-component wave functions depending on the three coordinates of the $i$-th 
electron,
we may infer that $H_{DC}$ acts on sixteen-component wave functions
which depend on the six coordinates of two electrons 
and have the anti-symmetric property due to the Pauli principle.

Mathematically, the operator should be written as
\begin{equation}
	H_{DC}=H(Z)\otimes I+I\otimes H(Z)+\frac{1}{| {\bf r}_{1}-{\bf r}_{2}|},
\end{equation} 
where 
\begin{equation}
	H(Z)=c\vec{\alpha}\cdot\vec{p}+mc^{2}\beta -\frac{Z}{|{\bf r}|}
\end{equation} 
is an operator in ${\cal H}=L^{2}({\bf R}^{3};\ {\bf C}^{4})$.
The domain of $H_{DC}$ is a subspace of the antisymmetric tensor product 
${\cal H}\otimes_{A} {\cal H}$.

In this paper we shall rigorously derive  a representation of $H_{DC}$ 
as a matrix operator
of order sixteen and give an answer to its spectral problems,
especially essential self-adjointness, continuous
(essential) spectrum and absence of eigenvalues. 
Our method was inspired by a simpler model operator 
(see (\ref{SM})).
\section{Two-electron problems}
As far as we know, there is no systematic derivation of relativistic systems
in the physics and quantum chemistry literature, so that
two-body relativistic systems with which we are concerned seem to be 
less familiar than nonrelativistic ones.
The first relativistic equation for two particles which
was extensively used in the past was the Breit equation; 
it is a differential equation for a relativistic wave function for two electrons,
interacting with each other and with an external electromagnetic field.
It is not fully Lorentz invariant and is only an approximation.
It reads (see \cite{BS})
\begin{equation}\label{Breit} 
	\left( E-H[{\bf 1}]-H[{\bf 2}]-\frac{e^{2}}{r_{12}}\right)U=
	-\frac{e^{2}}{2r_{12}}\left[ \vec{\alpha}_{1}\cdot \vec{\alpha}_{2}+
	\frac{(\vec{\alpha}_{1}\cdot {\bf r}_{12})
	(\vec{\alpha}_{2}\cdot {\bf r}_{12})}{r_{12}^{2}}\right]U,  
\end{equation} 
where ${\bf r}_{12}={\bf r}_{1}-{\bf r}_{2},\ r_{12}=|{\bf r}_{12}|$ and
\begin{equation} 
	H[{\bf j}]=-e\varphi ({\bf r}_{j})+\beta _{j}mc^{2}+\vec{\alpha}_{j}
	\cdot(c{\bf p}_{j}+e {\bf A}({\bf r}_{j}))	
\end{equation} 
is the Dirac Hamiltonian and the Dirac matrices $\vec{\alpha}_{j} $ and $\beta _{j}$
operate on the spinor of $U$ (for electron $j$).
The wave function $U$ depends on the positions ${\bf r}_{1}$ and ${\bf r}_{2}$
and has sixteen spinor components.

If we neglect the right hand side of the equation (\ref{Breit}) with ${\bf A}=0$,  
then we get the Dirac-Coulomb equation
\begin{equation}\label{DCeq} 
	\left( E-H_{D}[{\bf 1}]-H_{D}[{\bf 2}]-\frac{e^{2}}{r_{12}}\right)\Psi=0,
\end{equation} 
where $H_{D}[{\bf j}]$ is the usual Dirac operator acting on the $j$-th electron:
\begin{equation}
	H_{D}\Psi=\left(\begin{array}{cc}
	(m+V)I_{2}&\vec{\sigma }\cdot \vec{p}\\
	\vec{\sigma }\cdot \vec{p}&-(m-V)I_{2}
	\end{array}\right)
	\left(\begin{array}{c}
	\Psi^{\ell}\\
	\Psi^{s}
	\end{array}\right),\ \mathrm{with}\ V=-\frac{Z}{|{\bf r}|}.
\end{equation} 
Here, $\vec{\sigma }=(\sigma _{1},\sigma _{2},\sigma _{3})$ are the Pauli matrices
\[ 
	\sigma _{1}=\left(\begin{array}{cc}
	0&1\\
	1&0
	\end{array}\right),\
	\sigma _{2}=\left(\begin{array}{cc}
	0&-i\\
	i&0
	\end{array}\right),\
	\sigma _{3}=\left(\begin{array}{cc}
	1&0\\
	0&-1
	\end{array}\right)
\] 
and
\[ 
	\vec{\sigma }\cdot \vec{p}=\sum_{i=1}^{3}\sigma _{i}p_{i}.
\] 
The letters $\ell$ and $s$ refer to the large and the small part
of the wave function.

In the relativistic theory, we have to handle
a two-fold tensor product space of ${\bf C}^{4}$-valued functions
because the usual Dirac operator acts on four-vectors belonging
to $L^{2}({\bf R}^{3};\ {\bf C}^{4})$.
For ${\cal H}=L^{2}({\bf R}^{3};\ {\bf C}^{4})$, the two-fold tensor product 
${\cal H}^{\otimes 2}={\cal H}\otimes {\cal H}$ can be identified with
\[ 
	{\cal H}\otimes {\cal H}=\left\{ \psi({\bf 1},{\bf 2})
		=\;^{t}\! (\vec{\psi}_{\ell\ell},\vec{\psi}_{\ell s},
		\vec{\psi}_{s\ell},\vec{\psi}_{ss})
		\in L^{2}({\bf R}^{6};\ {\bf C}^{16})\ |\
		\vec{\psi}_{ij}\in L^{2}({\bf R}^{6};\ {\bf C}^{4})
	\right\} .
\]
Here we have identified $L^{2}({\bf R}^{6};\ {\bf C}^{4})$ with
$
	L^{2}({\bf R}^{6};\ {\bf C}^{2})\otimes {\bf C}^{2}
$
in the following way.

For $k=1,2$, let
\begin{align*}
	L^{2}({\bf R}^{3};\ {\bf C}^{4})&\ni\psi({\bf k})
	=\left(\begin{array}{c}
	\psi_{1}^{\ell}({\bf r}_{k})\\
	\psi_{2}^{\ell}({\bf r}_{k})\\
	\psi_{1}^{s}({\bf r}_{k})\\
	\psi_{2}^{s}({\bf r}_{k})
	\end{array}\right)
	=\Big(\psi_{1}^{\ell}e_{1}+\psi_{2}^{\ell}e_{2}\Big)\otimes f_{\ell}
	+
	\Big(\psi_{1}^{s}e_{1}+\psi_{2}^{s}e_{2}\Big)\otimes f_{s},\\
	\intertext{where}
	e_{1}=\; ^{t}\!(1,0),&\
	e_{2}=\; ^{t}\!(0,1),\ f_{\ell}=\; ^{t}\!(1,0),\
	f_{s}=\; ^{t}\!(0,1).
\end{align*}
In this notation, we see that any product function
\begin{equation}
	\psi({\bf 1})\otimes \psi({\bf 2})
	=\sum_{a,b\in \left\{ \ell,s\right\} }
	\sum_{i=1}^{2}\sum_{j=1}^{2} 
	\psi_{a,b,i,j}({\bf r}_{1},{\bf r}_{2})(e_{i}\otimes e_{j})\otimes 
	(f_{a}\otimes f_{b})
\end{equation}
satisfies
\begin{align*} 
	\vec{\psi}_{ab}({\bf r}_{1},{\bf r}_{2})
	=\left(\begin{array}{c}
	\psi_{a,b,1,1}({\bf r}_{1},{\bf r}_{2})\\
	\psi_{a,b,1,2}({\bf r}_{1},{\bf r}_{2})\\
	\psi_{a,b,2,1}({\bf r}_{1},{\bf r}_{2})\\
	\psi_{a,b,2,2}({\bf r}_{1},{\bf r}_{2})
	\end{array}\right)
	=\left(\begin{array}{c}
	\psi_{1}^{a}({\bf r}_{1})\\
	\psi_{2}^{a}({\bf r}_{1})\\
	\end{array}\right)\otimes \left(\begin{array}{c}
	\psi_{1}^{b}({\bf r}_{2})\\
	\psi_{2}^{b}({\bf r}_{2})
	\end{array}\right)
	=\left(\begin{array}{c}
	\psi_{1}^{a}({\bf r}_{1})\psi_{1}^{b}({\bf r}_{2})\\
	\psi_{1}^{a}({\bf r}_{1})\psi_{2}^{b}({\bf r}_{2})\\
	\psi_{2}^{a}({\bf r}_{1})\psi_{1}^{b}({\bf r}_{2})\\
	\psi_{2}^{a}({\bf r}_{1})\psi_{2}^{b}({\bf r}_{2})\\
	\end{array}\right)
\end{align*} 
for any $a,\ b\in \left\{ \ell,\ s\right\} $.

Now  we shall define two subspaces of 
${\cal H}^{\otimes 2}={\cal H}\otimes {\cal H}$,
the anti-symmetric space and symmetric space, denoted by 
${\cal H}_{A}^{2}={\cal H}\otimes _{A}{\cal H}$ and
${\cal H}_{S}^{2}={\cal H}\otimes _{S}{\cal H}$ , respectively.
\begin{definition}\label{AS20}
\begin{align*}
	{\cal H}_{A}^{2}=&\Big\{
		\psi({\bf r}_{1},{\bf r}_{2})
		=^{t}\! (\vec{\psi}_{11},\vec{\psi}_{12},
		\vec{\psi}_{21},\vec{\psi}_{22})
		\in L^{2}({\bf R}^{6};\ {\bf C}^{16})\ |\
		\vec{\psi}_{ij}\in L^{2}({\bf R}^{6};\ {\bf C}^{4}),\\
		&\psi_{ij}=Mat[\vec{\psi}_{ij}]\in L^{2}({\bf R}^{6};\ M(2,{\bf C})),
		\
	\psi_{k,k}({\bf r}_{2},{\bf r}_{1})
	=-^{t}\!\psi_{k,k}({\bf r}_{1},{\bf r}_{2}),\ k=1,2,\\
	&\psi_{1,2}({\bf r}_{2},{\bf r}_{1})
		=-^{t}\!\psi_{2,1}({\bf r}_{1},{\bf r}_{2})
	\Big\},
\end{align*}
where 
\begin{equation}
	Mat[\vec{a}]=\left(\begin{array}{cc}
	a_{1}&a_{3}\\
	a_{2}&a_{4}
	\end{array}\right)\ \mathrm{for}\ \vec{a}=\left(\begin{array}{c}
	a_{1}\\
	a_{2}\\
	a_{3}\\
	a_{4}
	\end{array}\right)\in {\bf C}^{4}
\end{equation} 
and $^{t}\!M$ stands for the transposed matrix of the $2\times 2$ matrix $M$.
\end{definition}
\begin{remark}
The definition \ref{AS20} coincides with the one in the literature 
on quantum chemistry (\cite{PBG}, \cite{PK}), where a slightly different notation from ours is adopted.
\begin{align}  \label{chem}
	\Psi({\bf 1},{\bf 2})=\left(\begin{array}{c}
	\vec{\Psi}^{\ell\ell}({\bf 1},{\bf 2})\\
	\vec{\Psi}^{\ell s}({\bf 1},{\bf 2})\\
	\vec{\Psi}^{s\ell}({\bf 1},{\bf 2})\\
	\vec{\Psi}^{ss}({\bf 1},{\bf 2})
	\end{array}\right)\
\mathrm{with}
\begin{array}{c}
	\vec{\Psi}^{\ell\ell}({\bf 1},{\bf 2})=-\vec{\Psi}^{\ell\ell}({\bf 2},{\bf 1})\\
	\vec{\Psi}^{\ell s}({\bf 1},{\bf 2})=-\vec{\Psi}^{s\ell}({\bf 2},{\bf 1})\\
	\vec{\Psi}^{ss}({\bf 1},{\bf 2})=-\vec{\Psi}^{ss}({\bf 2},{\bf 1}).
\end{array}
\end{align}
Moreover, it should be pointed out that the four components $\vec{\Psi}^{\ell\ell},\
\vec{\Psi}^{\ell s},\ \vec{\Psi}^{s\ell},\ \vec{\Psi}^{ss}$ in (\ref{chem})
are not functions from ${\bf R}^{3}\otimes {\bf R}^{3}$ to ${\bf C}^{4}$, but
they are functions from $({\bf R}^{3}\times \left\{\uparrow,\ \downarrow\right\})\otimes ({\bf R}^{3}\times \left\{\uparrow,\ \downarrow\right\}) $ to ${\bf C}$.
\end{remark}

In a similar way, we can define the symmetric tensor product space.
\begin{definition}\label{S2}
\begin{align*}
	{\cal H}_{S}^{2}=&\Big\{
		\psi({\bf r}_{1},{\bf r}_{2})=^{t}\! 
		(\vec{\psi}_{11},\vec{\psi}_{12},
		\vec{\psi}_{21},\vec{\psi}_{22})
		\in L^{2}({\bf R}^{6};\ {\bf C}^{16})\ |\
		\vec{\psi}_{ij}\in L^{2}({\bf R}^{6};\ {\bf C}^{4}),\\
		&\psi_{ij}=Mat[\vec{\psi}_{ij}]\in L^{2}({\bf R}^{6};\ M(2,{\bf C})),
		\
	\psi_{k,k}({\bf r}_{2},{\bf r}_{1})
	=^{t}\!\psi_{1,1}({\bf r}_{1},{\bf r}_{2}),\ k=1,2,
	\\
	&\psi_{1,2}({\bf r}_{2},{\bf r}_{1})=^{t}\!\psi_{2,1}({\bf r}_{1},{\bf r}_{2})
	\Big\}.
\end{align*}
\end{definition}
In the matrix formulation, we define the standard inner product of ${\cal H}_{A}^{2}$ or ${\cal H}_{S}^{2}$ by
\begin{equation}
	\langle \vec{F},\vec{G}\rangle =\sum_{i,j=1}^{2}\int_{{\bf R}^{6}}
	\mathrm{tr}\big(F_{ij}({\bf r}_{1},{\bf r}_{2})^{t}\! \overline{G_{ij}({\bf r}_{1},{\bf r}_{2})}\big)d{\bf r}_{1}d{\bf r}_{2}.
\end{equation} 
The notation $\langle \vec{F}_{ij},\vec{G}_{ij}\rangle $ will be used as
an abbreviation for the integral over the trace.
\section{The two-electron Dirac-Coulomb Hamiltonian}
As usual, the free and the one-particle Dirac operator are denoted by
\begin{equation}
	H_{0}=\vec{\alpha} \cdot \vec{p}+m\beta,\
	H=H_{0}+\frac{k}{|{\bf r}|}I_{4}.
\end{equation} 
Here $\vec{p}=\;^{t}\!(p_{1},p_{2},p_{3}),\ p_{j}=-i\partial _{{\bf r}_{j}}$;
$m$ is a nonnegative number and $k\in {\bf R}$.

Theorem \ref{VEC} prepares for a representation of
\begin{equation}
	H_{DC}=H_{D}\otimes I_{4}+I_{4}\otimes H_{D}
	+V_{0}({\bf r}_{1},{\bf r}_{2}),\
	V_{0}({\bf r}_{1},{\bf r}_{2})=k_{0}/|{\bf r}_{1}-{\bf r}_{2}|I_{16}
\end{equation} 
in the subspace ${\cal H}_{A}^{2}$ of ${\cal H}^{\otimes 2}={\cal H}\otimes {\cal H}$
which coincides with that used in literature on
quantum chemistry (\cite{KTA}, \cite{K}, \cite{PBG}).
\begin{lemma}Let 
\begin{equation}
	M_{j}=\left(\begin{array}{cc}
	B&A_{j}\\
	A_{j}&-B
	\end{array}\right)\in M(2,{\bf C}),\ j=1,2.
\end{equation} 
Then it holds that
\begin{equation}
	M_{1}\otimes I_{2}+I_{2}\otimes M_{2}
	=
	\left(\begin{array}{cccc}
	2B&A_{2}&A_{1}&0\\
	A_{2}&0&0&A_{1}\\
	A_{1}&0&0&A_{2}\\
	0&A_{1}&A_{2}&-2B
	\end{array}\right).
\end{equation} 
\end{lemma}
\Proof
In general, the Kronecker product of two matrices $X=(x_{ij})$ and $Y=(y_{ij})$ is defined by
\begin{equation}
	X\otimes Y=\Big(x_{ij}Y\Big),
\end{equation} 
so that $X\otimes Y$ is a matrix of size $mn\times k\ell$
if $X$ and $Y$ are $m\times n$ and $k\times \ell$ type, respectively.

\begin{align} 
	&M_{1}\otimes I_{2}
	=
	\left(\begin{array}{cccc}
	B&0&A_{1}&0\\
	0&B&0&A_{1}\\
	A_{1}&0&-B&0\\
	0&A_{1}&0&-B
	\end{array}\right),\
	&I_{2}\otimes M_{2}=
	\left(\begin{array}{cccc}
	B&A_{2}&0&0\\
	A_{2}&-B&0&0\\
	0&0&B&A_{2}\\
	0&0&A_{2}&-B
	\end{array}\right).
\end{align} 
\qed
\begin{theorem}\label{VEC}
Let ${\cal H}^{\otimes 2}\ni \psi=^{t}\! 
(\vec{\psi}_{11},\vec{\psi}_{12},\vec{\psi}_{21},\vec{\psi}_{22})
\in L^{2}({\bf R}^{6};\ {\bf C}^{4})^{4}$.
Then it holds that
\begin{equation}
	(H_{0}\otimes I_{4}+I_{4}\otimes H_{0})\psi=
	\left(\begin{array}{cccc}
	2mI_{4}&h_{2}&h_{1}&0\\
	h_{2}&0&0&h_{1}\\
	h_{1}&0&0&h_{2}\\
	0&h_{1}&h_{2}&-2mI_{4}
	\end{array}\right)
	\left(\begin{array}{c}
	\vec{\psi}_{11}\\
	\vec{\psi}_{12}\\
	\vec{\psi}_{21}\\
	\vec{\psi}_{22}
	\end{array}\right),
\end{equation}
where 
\begin{align} 
	h_{1}=(\vec{\sigma}\cdot \vec{p}_{1})\otimes I_{2},\
	h_{2}=I_{2}\otimes (\vec{\sigma}\cdot \vec{p}_{2}).
\end{align} 
\end{theorem}
\Proof
Denote $f_{\ell}$ and $f_{s}$ by ${\bf i}$ and ${\bf j}$,
respectively.
Let us consider two elements of ${\cal H}^{\otimes 2}$
\begin{equation}
	\Psi=\Psi_{1}({\bf r})\otimes {\bf i}
	+\Psi_{2}({\bf r})\otimes {\bf j},\
	\Psi^{\prime}=\Psi_{1}^{\prime}({\bf r})\otimes {\bf i}
	+\Psi_{2}^{\prime}({\bf r})\otimes {\bf j},
\end{equation}
where for $k=1,2$,
\begin{equation}
	\Psi_{k}({\bf r})=\left(\begin{array}{c}
	\psi_{k1}({\bf r})\\
	\psi_{k2}({\bf r})
	\end{array}\right),\
	\Psi_{k}^{\prime}({\bf r})=\left(\begin{array}{c}
	\psi_{k1}^{\prime}({\bf r})\\
	\psi_{k2}^{\prime}({\bf r})
	\end{array}\right).
\end{equation} 
We can regard $\Psi_{k}$ ($k=1,2$) as functions of ${\bf x}=({\bf r},\omega)$ as follows.
\begin{equation}
	\Psi_{k}({\bf r},\omega)=\psi_{k1}({\bf r})\chi_{+}(\omega)+
	\psi_{k2}({\bf r})\chi_{-}(\omega),
\end{equation} 
where $\chi_{\pm}$ are two orthonormal functions describing the spin
of electrons.

Then we have
\begin{equation}
	\Psi\otimes \Psi^{\prime}
	=\Psi_{1}\otimes \Psi_{1}^{\prime}\otimes ({\bf i}\otimes {\bf i})
	+\Psi_{1}\otimes \Psi_{2}^{\prime}\otimes ({\bf i}\otimes {\bf j})
	+\Psi_{2}\otimes \Psi_{1}^{\prime}\otimes ({\bf j}\otimes {\bf i})
	+\Psi_{2}\otimes \Psi_{2}^{\prime}\otimes ({\bf j}\otimes {\bf j}),
\end{equation} 
where
\begin{equation}
	\Psi_{k}\otimes \Psi_{\ell}^{\prime}=
	\left(\begin{array}{c}
	\psi_{k1}({\bf r}_{1})\psi_{\ell 1}^{\prime}({\bf r}_{2})\\
	\psi_{k1}({\bf r}_{1})\psi_{\ell 2}^{\prime}({\bf r}_{2})\\
	\psi_{k2}({\bf r}_{1})\psi_{\ell 1}^{\prime}({\bf r}_{2})\\
	\psi_{k2}({\bf r}_{1})\psi_{\ell 2}^{\prime}({\bf r}_{2})
	\end{array}\right)
	\cong
	\left(\begin{array}{cc}
	\psi_{k1}({\bf r}_{1})\psi_{\ell 1}^{\prime}({\bf r}_{2})&
	\psi_{k2}({\bf r}_{1})\psi_{\ell 1}^{\prime}({\bf r}_{2})\\
	\psi_{k1}({\bf r}_{1})\psi_{\ell 2}^{\prime}({\bf r}_{2})&
	\psi_{k2}({\bf r}_{1})\psi_{\ell 2}^{\prime}({\bf r}_{2})
	\end{array}\right).
\end{equation} 
We see that
\begin{align} 
	&(H_{0}\otimes I_{4})(\Psi\otimes \Psi^{\prime})
	=(H_{0}\Psi)\otimes \Psi^{\prime}\nonumber\\
	&=\big\{(\vec{\sigma} \cdot \vec{p} \Psi_{2}+mI_{2}\Psi_{1}){\bf i}
	+(\vec{\sigma} \cdot \vec{p} \Psi_{1}-mI_{2}\Psi_{2}){\bf j}\big\}
	\otimes \big(\Psi_{1}^{\prime}{\bf i}+\Psi_{2}^{\prime}{\bf j}\big)\nonumber\\
	&=[(\vec{\sigma} \cdot \vec{p} \Psi_{2}+mI_{2}\Psi_{1})
	\otimes \Psi_{1}^{\prime}]
	\otimes ({\bf i}\otimes {\bf i})
	+[(\vec{\sigma} \cdot \vec{p} \Psi_{2}+mI_{2}\Psi_{1})
	\otimes \Psi_{2}^{\prime}]
	\otimes ({\bf i}\otimes {\bf j})\nonumber\\
	&+[(\vec{\sigma} \cdot \vec{p} \Psi_{1}-mI_{2}\Psi_{2})
	\otimes \Psi_{1}^{\prime}]
	\otimes ({\bf j}\otimes {\bf i})
	+[(\vec{\sigma} \cdot \vec{p} \Psi_{1}-mI_{2}\Psi_{2})
	\otimes \Psi_{2}^{\prime}]
	\otimes ({\bf j}\otimes {\bf j}).
\end{align} 

Since the four vectors
$
	{\bf i}\otimes {\bf i},\ {\bf i}\otimes {\bf j},\
	{\bf j}\otimes {\bf i},\ {\bf j}\otimes {\bf j}
$
are linearly independent in ${\cal H}^{\otimes 2}$, we find
\begin{align} 
	&(H_{0}\otimes I_{4})(\Psi\otimes \Psi^{\prime})=\left(\begin{array}{c}
	(\vec{\sigma} \cdot \vec{p} \Psi_{2}+mI_{2}\Psi_{1})\otimes \Psi_{1}^{\prime}\\
	(\vec{\sigma} \cdot \vec{p} \Psi_{2}+mI_{2}\Psi_{1})\otimes \Psi_{2}^{\prime}\\
	(\vec{\sigma} \cdot \vec{p} \Psi_{1}-mI_{2}\Psi_{2})\otimes \Psi_{1}^{\prime}\\
	(\vec{\sigma} \cdot \vec{p} \Psi_{1}-mI_{2}\Psi_{2})\otimes \Psi_{2}^{\prime}
	\end{array}\right)\nonumber\\
	&=\left(\begin{array}{cccc}
	mI_{4}&0&(\vec{\sigma} \cdot \vec{p})\otimes I_{2}&0\\
	0&mI_{4}&0&(\vec{\sigma} \cdot \vec{p})\otimes I_{2}\\
	(\vec{\sigma} \cdot \vec{p})\otimes I_{2}&0&-mI_{4}&0\\
	0&(\vec{\sigma} \cdot \vec{p})\otimes I_{2}&0&-mI_{4}
	\end{array}\right)
	\left(\begin{array}{c}
	\Psi_{1}\otimes \Psi_{1}^{\prime}\\
	\Psi_{1}\otimes \Psi_{2}^{\prime}\\
	\Psi_{2}\otimes \Psi_{1}^{\prime}\\
	\Psi_{2}\otimes \Psi_{2}^{\prime}
	\end{array}\right).
\end{align}

Similarly, the identity
\begin{equation}
	(I_{4}\otimes H_{0})(\Psi\otimes \Psi^{\prime})
	=\Psi\otimes (H_{0}\Psi^{\prime})
	=\big(\Psi_{1}{\bf i}+\Psi_{2}{\bf j}\big)\otimes 
	\big\{(\vec{\sigma} \cdot \vec{p} \Psi_{2}^{\prime}+mI_{2}\Psi_{1}^{\prime}){\bf i}
	+(\vec{\sigma} \cdot \vec{p} \Psi_{1}^{\prime}-mI_{2}\Psi_{2}^{\prime}){\bf j}\big\}
\end{equation} 
implies
\begin{align} 
	&(I_{4}\otimes H_{0})(\Psi\otimes \Psi^{\prime})\nonumber\\
	&=\left(\begin{array}{cccc}
	mI_{4}&I_{2}\otimes(\vec{\sigma} \cdot \vec{p})&0&0\\
	I_{2}\otimes(\vec{\sigma} \cdot \vec{p})&-mI_{4}&0&0\\
	0&0&mI_{4}&I_{2}\otimes(\vec{\sigma} \cdot \vec{p})\\
	0&0&I_{2}\otimes(\vec{\sigma} \cdot \vec{p})&-mI_{4}
	\end{array}\right)
	\left(\begin{array}{c}
	\Psi_{1}\otimes \Psi_{1}^{\prime}\\
	\Psi_{1}\otimes \Psi_{2}^{\prime}\\
	\Psi_{2}\otimes \Psi_{1}^{\prime}\\
	\Psi_{2}\otimes \Psi_{2}^{\prime}
	\end{array}\right).
\end{align} 
\qed
\section{Main results}
The following spectral properties of the usual Dirac operator
\begin{align}
	 H=\vec{\alpha}\cdot \vec{p}+m\beta -k/|{\bf r}|I_{4}
\end{align}
with $m>0$ are well-known (see, e.g., \cite{T0}, \cite{W}).
\begin{enumerate}
\item 
$H$ is essentially selfadjoint on 
$C_{0}^{\infty }({\bf R}^{3})^{4}$ if $|k|\leq \sqrt{3}/2$.
\item
$\sigma _{ess}(H)={\bf R}\backslash (-m,m)$
if $|k|\leq \sqrt{3}/2$.
\item If $|k|\leq \sqrt{3}/2$, then
$H$ has no eigenvalues in ${\bf R}\backslash (-m,m)$ and 
there are countably many eigenvalues in $(-m,m)$ whose only 
accumulating points are $\pm m$.
\end{enumerate}
When the scalar potential $k/|{\bf r}|$ is replaced by
any symmetric matrix potential $V({\bf r})$ satisfying
$|V({\bf r})|\leq k/|{\bf r}|$, we have
\begin{enumerate}
\setcounter{enumi}{3}
\item
$
	\vec{\alpha }\cdot\vec{p}+m\beta +V({\bf r})
$
is essentially selfadjoint on $C_{0}^{\infty }({\bf R}^{3})^{4}$
if $|k|<1/2$.
\end{enumerate}

Now we shall state our main results for the Hamiltonian
\begin{equation}
	H_{DC}=H\otimes I_{4}+I_{4}\otimes H
	+\frac{k_{0}}{|{\bf r}_{1}-{\bf r}_{2}|}I_{16}.
\end{equation} 
These are just the first attempts to answer the spectral problem.
\begin{theorem}\label{ESA}Suppose that
$|k|<\sqrt{3}/2$. Then, for any nonzero real $k_{0}$, $H_{DC}$ is essentially selfadjoint on 
$[C_{0}^{\infty }({\bf R}^{3}\; ;\ {\bf C}^{4})]^{\otimes 2}\cap {\cal H}_{A}^{2}$.
\end{theorem}
\begin{remark}
The same conclusion is true if we replace the Coulomb potentials $k/|{\bf r}_{j}|$
by any symmetric matrix potentials $V({\bf r}_{j})$ satisfying
\begin{equation}
	|V({\bf r}_{j})|\leq k^{\prime}/|{\bf r}_{j}|,\ j=1,\ 2.
\end{equation} 
with $|k^{\prime}|<1/2$.
\end{remark}
The unique self-adjoint extension is denoted by the same symbol
$H_{DC}$ again. 
\begin{theorem} \label{ESS}
Suppose that $|k|<\sqrt{3}/2$. Then, for any nonzero real $k_{0}$, we have
\begin{equation}
	\sigma _{ess}(H_{DC}|_{{\cal H}_{A}^{2}})={\bf R}.
\end{equation} 
\end{theorem}
\begin{theorem}\label{NOE}
Suppose that $|k|<\sqrt{3}/2$. Then, for any nonzero real $k_{0}$, $H_{DC}|_{{\cal H}_{A}^{2}}$ has no eigenvalues.
\end{theorem}

We also consider the following simple model operator $\mathbb{H}$
\begin{equation}\label{SM}
	\mathbb{H}:=\alpha p_1+\alpha p_2+ 2m\beta + \frac{k_1}{|{\bf r}_{1}|}
	+\frac{k_2}{|{\bf r}_{2}|}
	+ \frac{k_0}{|{\bf r}_{1}-{\bf r}_{2}|}
\end{equation}
in $L^{2}({\bf R}^{6})^{4}$. It is easier to handle
than $H_{DC}$ because an orthogonal change of variables reduces it to
a Dirac operator in $L^{2}({\bf R}^{3})^{4}$ with a double-well potential.
\begin{theorem}\label{esa4} If $|k_j| < \sqrt{3}/2$ $(j=1,2)$, 
$\mathbb{H}$ on $C_0^\infty({\bf R}^6)^4$ is essentially self-adjoint for any 
$k_0 \in {\bf R}$.
\end{theorem}
Let $\mathbb{H}$ denote the unique self-adjoint extension again. 
\begin{theorem}\label{esssp4}
Suppose $k_0\neq 0$. Then the essential spectrum covers the whole line, that is, $\sigma_{{\rm  ess}}(\mathbb{H})= {\bf R}.$
\end{theorem}
\begin{theorem}\label{noev4} Let $k_0\neq 0$. Then
$\mathbb{H}$ has no eigenvalues, that is, $\sigma_{{\rm  p}}(\mathbb{H})= \emptyset$.    
\end{theorem} 
\section{The canonical form of $H_{DC}$ on the anti-symmetric space}
We return to the familiar notation ${\bf x}_{j}$ instead of ${\bf r}_{j}$.
We may represent $H_{DC}$ as follows.
Recall 
\begin{align} \label{smh}
	h_{1}=(\vec{\sigma}\cdot \vec{p}_{1})\otimes I_{2},\
	h_{2}=I_{2}\otimes (\vec{\sigma}\cdot \vec{p}_{2}).
\end{align} 
\begin{align} \label{DC}
	&H_{DC}\Psi=
	\left(\begin{array}{cccc}
	(2m+V)I_{4}&h_{2}&h_{1}&0\\
	h_{2}&VI_{4}&0&h_{1}\\
	h_{1}&0&VI_{4}&h_{2}\\
	0&h_{1}&h_{2}&-(2m-V)I_{4}
	\end{array}\right)
	\left(\begin{array}{c}
	\vec{\psi}_{11}({\bf x}_{1},{\bf x}_{2})\\
	\vec{\psi}_{12}({\bf x}_{1},{\bf x}_{2})\\
	\vec{\psi}_{21}({\bf x}_{1},{\bf x}_{2})\\
	\vec{\psi}_{22}({\bf x}_{1},{\bf x}_{2})
	\end{array}\right),\\
	&V({\bf x}_{1},{\bf x}_{2})
	=V({\bf x}_{1})+V({\bf x}_{2})+V_{0}({\bf x}_{1},{\bf x}_{2}).
\end{align}
\begin{proposition}\label{MAT}
Let $\Psi\in {\cal H}^{\otimes 2}=L^{2}({\bf R}^{6};\ M(2,{\bf C}))^{4}$.
Then
\begin{equation}
	H_{DC}\Psi=
	\left(\begin{array}{cccc}
	(2m+V)I_{4}&(h)_{2}&(h)_{1}&0\\
	(h)_{2}&VI_{4}&0&(h)_{1}\\
	(h)_{1}&0&VI_{4}&(h)_{2}\\
	0&(h)_{1}&(h)_{2}&-(2m-V)I_{4}
	\end{array}\right)
	\left(\begin{array}{c}
	\psi_{11}({\bf x}_{1},{\bf x}_{2})\\
	\psi_{12}({\bf x}_{1},{\bf x}_{2})\\
	\psi_{21}({\bf x}_{1},{\bf x}_{2})\\
	\psi_{22}({\bf x}_{1},{\bf x}_{2})
	\end{array}\right),
\end{equation}
where 
\begin{align} 
	&(h)_{1}\psi_{ij}({\bf x}_{1},{\bf x}_{2})
	=Mat\Big[\big((\vec{\sigma}\cdot \vec{p}_{1})\otimes I_{2}\big)
	 \vec{\psi}_{ij}({\bf x}_{1},{\bf x}_{2})\Big]
	=\vec{p}_{1}\psi_{ij}({\bf x}_{1},{\bf x}_{2})\cdot \; ^{t}\!\vec{\sigma} ,
	\label{mat1}\\
	&(h)_{2}\psi_{ij}({\bf x}_{1},{\bf x}_{2})
	=Mat\Big[\big((I_{2}\otimes \vec{\sigma}\cdot \vec{p}_{2})\big)\vec{\psi}_{ij}({\bf x}_{1},{\bf x}_{2})\Big]
	=\vec{p}_{2}\cdot \vec{\sigma} \psi_{ij}({\bf x}_{1},{\bf x}_{2}).\label{mat2}
\end{align} 
\end{proposition}
\Proof
Let
\begin{equation}
	\mathrm{vec}\left(\begin{array}{cc}
	a&b\\
	c&d
	\end{array}\right)=\left(\begin{array}{c}
	a\\
	c\\
	b\\
	d
	\end{array}\right).
\end{equation} 
Since
\begin{align} 
	&Mat\big[(A\otimes I_{2})\mathrm{vec} M\big]=M\; ^{t}\!A,\
	Mat[(I_{2}\otimes B)\mathrm{vec}M]=BM
\end{align} 
for any $2\times 2$ matrices $A$ and $M=(\Psi_{k}\otimes \Psi_{\ell}^{\prime})$,
we arrive at the conclusion.
\qed

\begin{proposition}
${\cal H}_{A}^{2}$ is an invariant subspace of ${\cal H}^{\otimes 2}={\cal H}\otimes {\cal H}$
with respect to the operator $H_{DC}$.
\end{proposition}
\Proof
If $\psi_{12}({\bf x}_{2},{\bf x}_{1})=-(\; ^{t}\! \psi_{21})({\bf x}_{1},{\bf x}_{2})$, then
$
	(\vec{p}_{2}\psi_{12})({\bf x}_{2},{\bf x}_{1})
	=-(\vec{p}_{1}\; ^{t}\! \psi_{21})({\bf x}_{1},{\bf x}_{2}).
$
We shall check the first component of $H_{DC}\Psi$.
In view of (\ref{mat1}) and \label{mat2} we find
\begin{align} 
	\big((h)_{2}\psi_{12}+&(h)_{1}\psi_{21}\big)({\bf x}_{2},{\bf x}_{1})
	=\vec{\sigma} \cdot \vec{p}_{2} \psi_{12}({\bf x}_{2},{\bf x}_{1})
	+\sum_{j=1}^{3}p_{1,j}\psi_{21}({\bf x}_{2},{\bf x}_{1})(^{t}\! \sigma _{j}) 
	\nonumber\\
	&=-\vec{\sigma} \cdot \vec{p}_{1} (^{t}\! \psi_{21})({\bf x}_{1},{\bf x}_{2})
	-\sum_{j=1}^{3}p_{2,j}(^{t}\! \psi_{12})({\bf x}_{1},{\bf x}_{2}) (\;^{t}\! \sigma _{j}) 
	\nonumber\\
	&=-^{t}\! \big((h)_{1}\psi_{21}\big)({\bf x}_{1},{\bf x}_{2})
	-^{t}\! \big((h)_{2}\psi_{12}\big)({\bf x}_{1},{\bf x}_{2})\nonumber\\
	&=-^{t}\!\big((h)_{2}\psi_{12}+(h)_{1}\psi_{21}\big)
	({\bf x}_{1},{\bf x}_{2}). 
\end{align} 
As for the second component, we see
\begin{align} 
	\big((h)_{2}\psi_{11}+&(h)_{1}\psi_{22}\big)({\bf x}_{2},{\bf x}_{1})
	=\vec{\sigma} \cdot \vec{p}_{2} \psi_{11}({\bf x}_{2},{\bf x}_{1})
	+\sum_{j=1}^{3}p_{1,j}\psi_{22}({\bf x}_{2},{\bf x}_{1})
	(^{t}\! \sigma _{j}) 
	\nonumber\\
	&=-\vec{\sigma} \cdot \vec{p}_{1} (^{t}\! \psi_{11})
	({\bf x}_{1},{\bf x}_{2})
	-\sum_{j=1}^{3}p_{2,j}(^{t}\! \psi_{22})({\bf x}_{1},{\bf x}_{2})
	(^{t}\! \sigma _{j}) 
	\nonumber\\
	&=-^{t}\! \big((h)_{1}\psi_{11}\big)({\bf x}_{1},{\bf x}_{2})
	-^{t}\! \big((h)_{2}\psi_{22}\big)({\bf x}_{1},{\bf x}_{2})\nonumber\\
	&=-^{t}\!\big((h)_{2}\psi_{11}+(h)_{1}\psi_{22}\big)
	({\bf x}_{1},{\bf x}_{2}) .
\end{align} 
As for the third and fourth components, similar computations yield
\begin{equation}
	\big((h)_{1}\psi_{11}+(h)_{2}\psi_{22}\big)({\bf x}_{2},{\bf x}_{1})
	=-^{t}\!\big((h)_{1}\psi_{11}+(h)_{2}\psi_{22}\big)
	({\bf x}_{1},{\bf x}_{2}) 
\end{equation} 
and
\begin{equation}
	\big((h)_{1}\psi_{12}+(h)_{2}\psi_{21}\big)({\bf x}_{2},{\bf x}_{1})
	=-^{t}\!\big((h)_{1}\psi_{12}+(h)_{2}\psi_{21}\big)
	({\bf x}_{1},{\bf x}_{2}).
\end{equation} 
\qed

The anti-symmetric property implies the following identities.
\begin{lemma}\label{INT}
Suppose that $F,G\in {\cal H}_{A}^{2}$. Then for any quadruple of indices
$i,\ j,\ k,\ \ell$,
\begin{align} 
	&\langle (\vec{\sigma }\cdot\vec{p}_{\tau}\otimes I_{2})\vec{F}_{ij},\vec{G}_{k\ell}\rangle 
	=\langle (I_{2}\otimes \vec{\sigma }\cdot\vec{p}_{\tau})\vec{F}_{ji},\vec{G}_{\ell k}\rangle ,\ for\  \tau =1,2,\\
	&\langle V\vec{F}_{ij},\vec{G}_{ij}\rangle =\langle V\vec{F}_{ji},\vec{G}_{ji}\rangle 
\end{align} 
if $V$ is a diagonal matrix satisfying
\begin{equation}
	V({\bf x}_{1},{\bf x}_{2})=V({\bf x}_{2},{\bf x}_{1}).
\end{equation}
\end{lemma}
\Proof For simplicity, we only consider the case when $\tau=2$.
Both $\psi\in {\cal H}_{S}^{2}$ and $\varphi \in {\cal H}_{A}^{2}$ satisfy
\begin{align} 
\begin{split}
	&\partial_{{\bf x}_{2}}\vec{\psi}_{ij}({\bf x}_{1},{\bf x}_{2})
	=\mathrm{vec}[\partial_{{\bf x}_{1}}\;^{t}\!\psi_{ji}]({\bf x}_{2},{\bf x}_{1}),\
\\	&\partial_{{\bf x}_{2}}\vec{\varphi} _{ij}({\bf x}_{1},{\bf x}_{2})
	=-\mathrm{vec}[\partial_{{\bf x}_{1}}\;^{t}\!\varphi_{ji}]({\bf x}_{2},{\bf x}_{1}).
\end{split}\end{align} 
for any $i,j=1,2$.

Hence for any $\psi,\ \varphi \in {\cal H}_{S}^{2}$,
\begin{align} 
	&\langle I_{2}\otimes \vec{\sigma} \cdot\nabla_{{\bf x}_{2}}\vec{\psi}_{ij},
	\vec{\varphi}_{k\ell} \rangle \nonumber\\
	&=\int_{{\bf R}^{6}}tr\big((\vec{\sigma}\cdot\nabla_{{\bf x}_{2}}\psi_{ij})
	({\bf x}_{1},{\bf x}_{2})
	\;^{t}\!\overline{\varphi}_{k\ell} ({\bf x}_{1},{\bf x}_{2})\big)d{\bf x}_{1}d{\bf x}_{2}
	\nonumber\\
	&= \int_{{\bf R}^{6}}tr\big(\vec{\sigma}\cdot\nabla_{{\bf x}_{1}}
	\;^{t}\!\psi_{ji}({\bf x}_{2},{\bf x}_{1})
	\;\overline{\varphi}_{\ell k} ({\bf x}_{2},{\bf x}_{1})
	\big)d{\bf x}_{1}d{\bf x}_{2}\nonumber\\
	&=\int_{{\bf R}^{6}}tr\big((\nabla_{{\bf x}_{1}}\psi_{ji}({\bf x}_{2},{\bf x}_{1})
	\cdot^{t}\!\sigma)
	\;^{t}\!\overline{\varphi}_{\ell k} ({\bf x}_{2},{\bf x}_{1})
	\big)d{\bf x}_{2}d{\bf x}_{1}\nonumber\\
	&=\langle \vec{\sigma}\cdot \nabla_{{\bf x}_{2}}\otimes I_{2}\vec{\psi}_{ji},
	\vec{\varphi}_{\ell k} \rangle .
\end{align} 
Here, we have used
\begin{equation}
	tr(AB)=tr(^{t}\!(AB))=tr(^{t}\! B\;^{t}\! A)=tr(^{t}\! A\;^{t}\! B).
\end{equation} 
Similarly, we see that
for any $\psi,\ \varphi \in {\cal H}_{A}^{2}$,
\begin{align} 
	\langle I_{2}\otimes \vec{\sigma} \cdot\nabla_{{\bf x}_{2}}\vec{\psi}_{ij},
	\vec{\varphi}_{k\ell} \rangle 
	=\langle \vec{\sigma}\cdot \nabla_{{\bf x}_{2}}\otimes I_{2}\vec{\psi}_{ji},
	\vec{\varphi}_{\ell k}\rangle .
\end{align} 

\qed
\begin{theorem}\label{FORM}
Let
\begin{equation}
	H_{DC}^{+}\Psi=
	\left(\begin{array}{cccc}
	V+2m&h_{12}&0&0\\
	h_{12}&V&0&0\\
	0&0&V&h_{12}\\
	0&0&h_{12}&V-2m
	\end{array}\right)
	\left(\begin{array}{c}
	\vec{\psi}_{11}({\bf x}_{1},{\bf x}_{2})\\
	\vec{\psi}_{12}({\bf x}_{1},{\bf x}_{2})\\
	\vec{\psi}_{21}({\bf x}_{1},{\bf x}_{2})\\
	\vec{\psi}_{22}({\bf x}_{1},{\bf x}_{2})
	\end{array}\right)
\end{equation} 
and
\begin{equation}
	H_{DC}^{-}\Psi=
	\left(\begin{array}{cccc}
	V+2m&0&h_{21}&0\\
	0   &V&0     &h_{21}\\
	h_{21}&0&V&0\\
	0&h_{21}&0&V-2m
	\end{array}\right)
	\left(\begin{array}{c}
	\vec{\psi}_{11}({\bf x}_{1},{\bf x}_{2})\\
	\vec{\psi}_{12}({\bf x}_{1},{\bf x}_{2})\\
	\vec{\psi}_{21}({\bf x}_{1},{\bf x}_{2})\\
	\vec{\psi}_{22}({\bf x}_{1},{\bf x}_{2})
	\end{array}\right),
\end{equation} 
where
\begin{align*} 
	&h_{12}=I_{2}\otimes (\vec{\sigma}\cdot \vec{p}_{1})+
              I_{2}\otimes (\vec{\sigma}\cdot \vec{p}_{2}).\\
&h_{21}=(\vec{\sigma}\cdot \vec{p}_{1})\otimes I_{2}+
(\vec{\sigma}\cdot \vec{p}_{2})\otimes I_{2}.
\end{align*} 
Then, 
If $\Psi,\Phi\in {\cal H}_{A}^{2}\cap C_{0}^{\infty }$, then
\begin{equation}
	\langle H_{DC}\Psi,\Phi\rangle =\langle H_{DC}^{+}\Psi,\Phi\rangle ,\
	\langle H_{DC}\Psi,\Phi\rangle =\langle H_{DC}^{-}\Psi,\Phi\rangle .
\end{equation} 
\end{theorem}

{\bf Proof of Theorem \ref{FORM}}
We need to consider both (\ref{smh}) and their variants as follows.
\begin{align*} 
	\check{h}_{1}=I_{2}\otimes (\vec{\sigma}\cdot \vec{p}_{1}),\
         \check{h}_{2}=(\vec{\sigma}\cdot \vec{p}_{2})\otimes I_{2}.
\end{align*} 
By virtue of Lemma \ref{INT}, we have
\begin{align} 
	&\langle h_{1}\vec{\psi}_{21},\vec{\phi}_{11}\rangle 
	=\langle \check{h}_{1}\vec{\psi}_{12},\vec{\phi}_{11}\rangle ,\
	\langle h_{1}\vec{\psi}_{12},\vec{\phi}_{22}\rangle 
	=\langle \check{h}_{1}\vec{\psi}_{21},\vec{\phi}_{22}\rangle ,
\\
	&\langle h_{1}\vec{\psi}_{11},\vec{\phi}_{21}\rangle 
	=\langle \check{h}_{1}\vec{\psi}_{11},\vec{\phi}_{12}\rangle ,\
	\langle h_{1}\vec{\psi}_{22},\vec{\phi}_{12}\rangle 
	=\langle \check{h}_{1}\vec{\psi}_{22},\vec{\phi}_{21}\rangle ,
\end{align} 
which implies
\begin{equation}
	\langle H_{DC}\Psi,\Phi\rangle =\langle H_{DC}^{+}\Psi,\Phi\rangle .
\end{equation}
Similarly,
\begin{align*} 
	&\langle h_{2}\vec{\psi}_{12},\vec{\phi}_{11}\rangle 
	=\langle \check{h}_{2}\vec{\psi}_{21},\vec{\phi}_{11}\rangle ,\
	\langle h_{2}\vec{\psi}_{21},\vec{\phi}_{22}\rangle 
	=\langle \check{h}_{2}\vec{\psi}_{12},\vec{\phi}_{22}\rangle ,
\\
	&\langle h_{2}\vec{\psi}_{11},\vec{\phi}_{12}\rangle 
	=\langle \check{h}_{2}\vec{\psi}_{11},\vec{\phi}_{2}\rangle ,\
	\langle h_{2}\vec{\psi}_{22},\vec{\phi}_{21}\rangle 
	=\langle \check{h}_{2}\vec{\psi}_{22},\vec{\phi}_{1,2}\rangle ,
\end{align*} 
which implies
\begin{equation} 
	\langle H_{DC}\Psi,\Phi\rangle =\langle H_{DC}^{-}\Psi,\Phi\rangle .
\end{equation}
\qed
\section{Proof of the spectral properties of $H_{DC}$}
\subsection{Essential selfadjointness}
{\bf Proof of Theorem \ref{ESA}}:
In view of Theorem \ref{FORM}, we can identify $H_{DC}$ with $H_{DC}^{+}$.
\begin{lemma}
Let $\varphi (t)\in C^{\infty}({\bf R})\cap L^{\infty }({\bf R}) $. Then
\begin{equation}
	[H_{DC}^{+},\varphi (|{\bf x}_{1}-{\bf x}_{2}|)]=0.
\end{equation} 
\end{lemma}

Let $\chi\in C_{0}^{\infty}({\bf R})$ satisfy that $0\leq \chi \leq 1$,  
\begin{equation}
	\chi(t)=\Big\{ 
	\begin{array}{ll}1, & t\geq 2, \\ 0, & t \leq 1, \end{array}
\end{equation}
and $B_{n}$ be a multiplication operator  defined by $\displaystyle 
	B_{n}=\chi(n|{\bf x}_{1}-{\bf x}_{2}|)
$ and
\begin{equation}
	V_{n}=B_{n}VB_{n}.
\end{equation} 
\begin{lemma}\label{APP}
$H_{DC,n}=H_{DC,0}+V_{n}$ is essentially selfadjoint on
$[C_{0}^{\infty }({\bf R}^{3};\ {\bf C}^{4})]^{\otimes 2}\cap {\cal H}_{A}^{2}$.
\end{lemma}
\begin{lemma}\label{THA}
\begin{equation}
 \lim_{n \rightarrow \infty} (H^* B_n\psi, \psi-B_n\psi)=0 
\end{equation} 
for any $\psi \in D(H^*)$.
\end{lemma}
Thanks to Theorem 5.2 of Thaller \cite{T}, from Lemma \ref{APP} and Lemma \ref{THA},
it follows that $H_{DC}$ is also essentially selfadjoint on
$[C_{0}^{\infty }({\bf R}^{3};\ {\bf C}^{4})]^{\otimes 2}\cap {\cal H}_{A}^{2}$.

Now we are in a position to prove Lemma \ref{APP}.
We consider an orthogonal transformation $S$ in $M(8,{\bf C})$
\begin{equation}
	S=\frac{1}{\sqrt{2}}\left(\begin{array}{cc}
	I_{4}&I_{4}\\
	I_{4}&-I_{4}
	\end{array}\right).
\end{equation} 
Then
\begin{align} 
	&S^{-1}\left(\begin{array}{cc}
	a&0\\
	0&b
	\end{array}\right)S=\frac{1}{2}\left(\begin{array}{cc}
	a+b&a-b\\
	a-b&a+b
	\end{array}\right),\\
	&S^{-1}\left(\begin{array}{cc}
	0&a\\
	a&0
	\end{array}\right)S=\left(\begin{array}{cc}
	a&0\\
	0&-a
	\end{array}\right).
\end{align} 

With the orthogonal transformation
\begin{equation}
	T=S\oplus S
\end{equation} 
we therefore have
\begin{align} \label{cano}
	&T^{-1}H_{DC,n}T\nonumber\\
	&=
	\left(\begin{array}{cccc}
	h_{12}+V_{n}+m&m&0&0\\
	m&-h_{12}+V_{n}+m&0&0\\
	0&0&h_{12}+V_{n}-m&-m\\
	0&0&-m&-h_{12}+V_{n}-m
	\end{array}\right)\\
	&=\left(\begin{array}{cccc}
	H_{00}&0&0&0\\
	0&-H_{00}&0&0\\
	0&0&H_{00}&0\\
	0&0&0&-H_{00}
	\end{array}\right)+m\left(\begin{array}{cccc}
	I_{4}&I_{4}&0&0\\
	I_{4}&I_{4}&0&0\\
	0&0&-I_{4}&-I_{4}\\
	0&0&-I_{4}&-I_{4}
	\end{array}\right)+V_{n}I_{16},
\end{align} 
where
\begin{align} 
	&H_{00}=h_{12}=I_{2}\otimes \big[\vec{\sigma}\cdot (\vec{p}_{1}+\vec{p}_{2})\big].
\end{align} 
Introducing a change of coordinates
\begin{align*} 
	\frac{1}{\sqrt{2}}({\bf x}_{1}+{\bf x}_{2})={\bf y}_{1},\quad
        \frac{1}{\sqrt{2}}({\bf x}_{1}-{\bf x}_{2})={\bf y}_{2},
\end{align*}
we find
\begin{align} 
	&H_{00}=\sqrt{2}[I_{2}\otimes (\vec{\sigma}\cdot \vec{p}_{{\bf y}_{1}})],\\
	&V=\frac{\sqrt{2}k}{|{\bf y}_{1}+{\bf y}_{2}|}+\frac{\sqrt{2}k}{|{\bf y}_{1}-{\bf y}_{2}|}
	+\frac{k_{0}}{\sqrt{2}|{\bf y}_{2}|}.
\end{align} 
Let ${\bf y}_{1}=(\eta_{1},\eta_{2},\eta_{3})$, $\vec{p}=(p_{1},p_{2},p_{3})$
with $p_{j}=-i\partial _{\eta_{j}}$ and
$\vec{\sigma} \cdot \vec{p}=\sum_{j=1}^{3}\sigma _{j}p_{j}$,
then it holds that
\begin{align} 
	H_{00}^{2}=2\big(I_{2}\otimes (\vec{\sigma} \cdot \vec{p})\big)^{2}=
	2|\vec{p}|^{2}I_{4}.
\end{align} 

We note that $H_{00}$ satisfies similar estimates to the ones for the Dirac operator
$\vec{\alpha }\cdot \vec{p}\;$ (see \cite{Sc0}, \cite{Sc}).
\begin{lemma}\label{WIN}
\begin{equation}
	\||{\bf y}_{1}|^{1/2}H_{00}u\|_{L^{2}({\bf R}^{6})^{4}}\geq
	\||{\bf y}_{1}|^{-1/2}u\|_{L^{2}({\bf R}^{6})^{4}}
\end{equation} 
provided that $u$ and $\nabla_{{\bf y}_{1}}u\in L^{2}({\bf R}^{6};\ {\bf C}^{4})$.
\end{lemma}
\begin{lemma}\label{HA}
Let $Q\in L^{\infty }_{\mathrm{loc}}\big({\bf R}^{3}\backslash\{0\}\big)^{4\times 4}$ be 
an Hermitian matrix which commutes with
$I_{2}\otimes (\vec{\sigma} \cdot \vec{{\bf y}}_{1})$ ($\vec{{\bf y}}_{1}\in {\bf R}^{3}$),
and which satisfies
\begin{equation}
	|Q({\bf y}_{1})|\leq 
	\sqrt{2}\mu |{\bf y}_{1}|^{-1},\ 0\leq \mu<\sqrt{3}/2.
\end{equation} 
Then for all $v\in H^{1}({\bf R}^{6};\ {\bf C}^{4})$ we have
\begin{equation}
	\|\sqrt{2}|{\bf y}_{1}|^{-1}v\|_{L^{2}({\bf R}^{6})^{4}}
	\leq \frac{1}{1-a}\|(H_{00}+Q({\bf y}_{1}))v\|_{L^{2}({\bf R}^{6})^{4}},\
	a=\sqrt{\mu^{2}+1/4}.
\end{equation} 
\end{lemma}
These inequalities can be proved by introducing polar coordinates in the variable ${\bf y}_{1}$
just as in the case of the usual Dirac operator (Schmincke [8,9]).

Without loss of generality, we may assume that $m=0$.
\begin{proposition}\label{ESSn}
Suppose that $|k|<\sqrt{3}/2$ and $k_{0}\neq 0$. Then,
$H_{n}:=H_{00}+V_{n}$ is essentially self-adjoint on $C_{0}^{\infty }({\bf R}^{6};\ {\bf C}^{4})$.
\end{proposition}
\Proof
It suffices to show that
\begin{equation}\label{Su} 
	\overline{(H_{n}\pm iI_{4})(C_{0}^{\infty }({\bf R}^{6};\ {\bf C}^{4}))}
	=L^{2}({\bf R}^{6};\ {\bf C}^{4}).
\end{equation} 

Due to the cutoff function,
the scalar potential $V_{n}$ has singularities which do not coincide, so that
we are able to employ a technique similar to the one developed by Vogelsang (\cite{V})
in view of the following lemma which can be shown to be based on Lemma \ref{WIN}.
\qed
\noindent
{\bf Proof of Lemma \ref{APP}}
From Proposition \ref{ESSn}, it follows that
\begin{equation}
	\overline{(H_{DC,n}\pm  iI_{4})
	\big([C_{0}^{\infty }({\bf R}^{6};\ {\bf C}^{4})]^{\otimes 2}\big)}
	=L^{2}({\bf R}^{6};\ {\bf C}^{16}).
\end{equation} 
Since
$[C_{0}^{\infty }({\bf R}^{3};\ {\bf C}^{4})]^{\otimes 2}\cap {\cal H}_{A}^{2}$ is 
an invariant subspace of $H_{DC,n}$, we can conclude that
$H_{DC,n}$ is essentially selfadjoint in ${\cal H}^{\otimes 2}_{A}$.
\qed
\subsection{  Essential spectrum}
{\bf Proof of Theorem \ref{ESS}}
Let $\lambda>m,\ \mu>m$. 
Take $\xi,\ \eta\in {\bf R}^{3}\backslash\{0\}$ such that
\begin{equation}
	|\xi|^{2}+m^{2}=\lambda^{2},\
	|\eta|^{2}+m^{2}=\mu^{2},\ \xi\cdot \eta=0.
\end{equation} 
For each $\xi\in {\bf R}^{3}$, let $u,v\in {\bf C}^{4}\backslash\{0\}$ be 
normalized eigenvectors to the equations
\begin{align} 
	&(\alpha \cdot\xi+m\beta )u=\lambda u,\
	|u|_{{\bf C}^{4}}=1,
\\
	&(\alpha \cdot\xi+m\beta )v=-\mu v,\
	|v|_{{\bf C}^{4}}=1,
\end{align} 
respectively.
Define two functions by
\begin{equation}
	u_{n}({\bf x})=\chi_{n}({\bf x})e^{i{\bf x}\cdot\xi}u(\xi),\
	v_{n}({\bf x})=\chi_{n}({\bf x})e^{i{\bf x}\cdot\eta}v(\xi),
\end{equation} 
where $\chi_{n}\in C_{0}^{\infty }({\bf R}^{3})$ is a nonnegative function
such that
\begin{equation}
	\int_{{\bf R}^{3}}|\chi_{n}({\bf x})|^{2}d{\bf x}=1,\
	\mathrm{supp}\chi_{n}\subset \{{\bf x} \in {\bf R}^{3} \ |\ n <|{\bf x}|< 2n \}.
\end{equation} 
Then $\{u_{n}({\bf x})\}$ and $\{v_{n}({\bf x})\}$ are two sets of
singular sequences in $L^{2}({\bf R}^{3})^{4}$
such that 
\begin{align}
	& \|u_{n}\|_{L^{2}}=\|v_{n}\|_{L^{2}}=1,\ u_{n}\rightharpoonup  0,\ 
	v_{n}\rightharpoonup  0,\ weakly,\\
	& [H_{0}+m\beta -\lambda]u_{n} \longrightarrow 0 \ \ (n\rightarrow\infty), 
	\quad {\rm supp}u_{n} \subset \{{\bf x} \in {\bf R}^{3} \ |\ n <|{\bf x}|< 2n \} \\
	& [H_{0}+m\beta +\mu]v_{n} \longrightarrow 0 \quad (n\rightarrow\infty),
	\quad {\rm supp}v_{n} \subset \{x\in {\bf R}^3\ | \ n <|{\bf x}|< 2n\}.
\end{align}
Now we shall construct a singular sequence in
${\cal H}_{A}^{2}$.
\begin{lemma}\label{ASS}
\begin{equation}
	w_{n}:=\frac{1}{\sqrt{2}}
	\{u_{n^{2}}\otimes v_{n}-v_{n}\otimes u_{n^{2}}\}\in {\cal H}_{A}^{2},
\end{equation}
\end{lemma}
\Proof
Let us consider $F=f_{1}\otimes {\bf i}+f_{2}\otimes {\bf j},\
 G=g_{1}\otimes {\bf i}+g_{2}\otimes {\bf j}
\in L^{2}({\bf R}^{3})^{4}$,
where $f_{j},\ g_{j}\in L^{2}({\bf R}^{3})^{2}$.
Recall the definition of $u\otimes v$:
\begin{equation}
	F\otimes G=
	(f_{1}\otimes g_{1})\otimes ({\bf i}\otimes {\bf i})
	+(f_{1}\otimes g_{2})\otimes ({\bf i}\otimes {\bf j})
	+(f_{2}\otimes g_{1})\otimes ({\bf j}\otimes {\bf i})
	+(f_{2}\otimes g_{2})\otimes ({\bf j}\otimes {\bf j}),
\end{equation} 
where $f_{j}\otimes g_{k}$ ($j,k=1,2$) are defined as
\begin{equation}
	f_{j}\otimes g_{k}=\left(\begin{array}{c}
	u_{1,j}\otimes v_{1,k}\\
	u_{1,j}\otimes v_{2,k}\\
	u_{2,j}\otimes v_{1,k}\\
	u_{2,j}\otimes v_{2,k}
	\end{array}\right)
	=
	\left(\begin{array}{c}
	u_{1,j}({\bf x}_{1}) v_{1,k}({\bf x}_{2})\\
	u_{1,j}({\bf x}_{1}) v_{2,k}({\bf x}_{2})\\
	u_{2,j}({\bf x}_{1}) v_{1,k}({\bf x}_{2})\\
	u_{2,j}({\bf x}_{1}) v_{2,k}({\bf x}_{2})
	\end{array}\right).
\end{equation}
for $\displaystyle  f_{j}=\left(\begin{array}{c}
	u_{1,j}\\
	u_{2,j}
	\end{array}\right)
	$
and $\displaystyle g_{k}=\left(\begin{array}{c}
	v_{1,k}\\
	v_{2,k}
	\end{array}\right)$.

It follows that
\begin{align}
	F\otimes G-G\otimes F
	&=(f_{1}\otimes g_{1}-g_{1}\otimes f_{1})\otimes ({\bf i}\otimes {\bf i})
	+(f_{1}\otimes g_{2}-g_{1}\otimes f_{2})\otimes ({\bf i}\otimes {\bf j})
	\nonumber\\
	&+(f_{2}\otimes g_{1}-g_{2}\otimes f_{1})\otimes ({\bf j}\otimes {\bf i})
	+(f_{2}\otimes g_{2}-g_{2}\otimes f_{2})\otimes ({\bf j}\otimes {\bf j})
	\nonumber\\
	&\in {\cal H}_{A}^{2}.
\end{align} 
\qed

From Lemma \ref{ASS}, it follows that
$\{w_{n}\}$ becomes a singular sequence of $H_{DC}-\lambda +\mu$.
In fact, it satisfies that as $n\to\infty $
\begin{align}
& (H_{DC}-(\lambda -\mu))w_{n}= 
\frac{k}{|{\bf x}_{1}|}w_{n}+ \frac{k}{|{\bf x}_{2}|}w_{n}+\frac{k_{0}}{|{\bf x}_{1}-{\bf x}_{2}|}w_{n}
\ \longrightarrow 0,\\
& w_{n}\rightharpoonup  0,\
\|w_{n}\|_{L^{2}}=1.
\end{align} 
The last assertion $\|w_{n}\|_{L^{2}}=1$ follows from the fact that the eigenvectors
$u(\xi),\ v(\xi)$ corresponding to the different eigenvalues
are orthogonal to each other.
Since
\[ 
	(m,\infty )+(-\infty ,-m)=(-\infty ,\infty ),
\]
we can conclude that $\sigma_{{\rm ess}}(H_{DC})={\bf R}$. 
\qed
\subsection{Absence of eigenvalues}
{\bf Proof of Theorem \ref{NOE}}
We first consider the case where $m>0$.

Let $u\in {\cal H}_{A}^{2}$ be a solution to
\begin{align} \label{EV}
	H_{DC}u=\lambda u.
\end{align} 
Then we see that for any function 
$\varphi({\bf y}_{1},{\bf y}_{2})
\in C_{0}^{\infty }({\bf R}^{6};\ {\bf C}^{16})\cap {\cal H}_{A}^{2}$, 
\begin{equation}\label{} 
	\langle H_{DC}u,\varphi \rangle _{{\cal H}_{A}^{2}}
	=\langle H_{DC}^{+}u,\varphi \rangle _{{\cal H}_{A}^{2}}.
\end{equation} 
Let 
\[ 
	(Fu)({\bf y})=u\big(({\bf y}_{1}+{\bf y}_{2})/\sqrt{2},
	({\bf y}_{1}-{\bf y}_{2})/\sqrt{2}\big).
\] 
Then we see that for $v\in C_{0}^{\infty }({\bf R}^{6};\ {\bf C}^{16})\cap {\cal H}_{A}^{2}$,
\begin{align}
	&\langle H_{DC}Fu,Fv\rangle _{{\cal H}_{A}^{2}}=\langle H_{DC}^{+}Fu,Fv\rangle _{{\cal H}_{A}^{2}},
\end{align} 
so that we can identify $H_{DC}$ with $H_{DC}^{+}$ in ${\cal H}_{A}^{2}$.
In view of (\ref{EV}), we have
\begin{align} \label{REV}
	H_{{\bf y}_{2}}(Fu)(\cdot,{\bf y}_{2})
	=\Big(\lambda -\frac{k_{0}}{\sqrt{2}|{\bf y}_{2}|}\Big) (Fu)(\cdot,{\bf y}_{2}),
\end{align} 
where $H_{{\bf y}_{2}}$ in $L^{2}({\bf R}^{3})^{16}$ is an operator with a parameter ${\bf y}_{2}$
defined by
\begin{align} 
	H_{{\bf y}_{2}}=\left(\begin{array}{cccc}
	H_{00}+mI_{4}&mI_{4}&0&0\\
	mI_{4}&-H_{00}+mI_{4}&0&0\\
	0&0&H_{00}-mI_{4}&-mI_{4}\\
	0&0&-mI_{4}&-H_{00}-mI_{4}
	\end{array}\right)+V_{{\bf y}_{2}}I_{16},
\end{align} 
with
\begin{align} 
	V_{{\bf y}_{2}}=\frac{\sqrt{2}k}{|{\bf y}_{1}+{\bf y}_{2}|}
	+\frac{\sqrt{2}k}{|{\bf y}_{1}-{\bf y}_{2}|}.
\end{align} 
Let
\begin{align} 
	H_{++}=\left(\begin{array}{cc}
	H_{00}+mI_{4}&mI_{4}\\
	mI_{4}&-H_{00}+mI_{4}\\
	\end{array}\right),\
	H_{--}=\left(\begin{array}{cc}
	H_{00}-mI_{4}&-mI_{4}\\
	-mI_{4}&-H_{00}-mI_{4}
	\end{array}\right).
\end{align} 
Then,
\begin{align} \label{square}
	(H_{++}-mI_{8})^{2}=(|p|^{2}+m^{2})I_{8},\
	(H_{--}+mI_{8})^{2}=(|p|^{2}+m^{2})I_{8}
\end{align}
with $p=-i\nabla_{{\bf y}_{1}}$. 
Since
\begin{align} 
	&\sigma _{ess}(H_{++}-mI_{8})=(-\infty ,-m]\cup [m,\infty ),\\
	&\sigma _{ess}(H_{--}+mI_{8})=(-\infty ,-m]\cup [m,\infty )
\end{align} 
and $\lim_{|{\bf y}_{1}|\to \infty }V_{{\bf y}_{2}}=0$, it holds that
\begin{align} 
	&\sigma _{ess}(H_{++}+V_{{\bf y}_{2}}I_{4})= (-\infty ,0]\cup [2m,\infty ),\\
	 &\sigma _{ess}(H_{--}+V_{{\bf y}_{2}}I_{4})=(-\infty ,-2m]\cup [0,\infty ).
\end{align}
On the other hand, because of Theorem \ref{ESA}, 
a technique of Weidmann (\cite{W} Theorem 10.38) and Kalf \cite{Kalf}
enables us to prove that
$H_{\pm\pm}+V_{{\bf y}_{2}}I_{4}$ have no eigenvalues in $(-\infty ,0]\cup [2m,\infty )$ and
$(-\infty ,-2m]\cup [0,\infty )$, respectively. 
Therefore, we see that  every eigenvalue of $H_{{\bf y}_{2}}$
is in $(-2m,0)\cup (0,2m)$ and its multiplicity
is finite if it exists.

For any fixed ${\bf y}_{2}\in{\bf R}^{3}\backslash\{0\}$ and $\kappa>0$, define an 
operator in  $L^{2}({\bf R}^{3})^{16}$
\[ 
	H(\kappa)=H_{\kappa{\bf y}_{2}}.
\]
The family of operators $\left\{  H(\kappa)\right\}_{\kappa>0}$
forms an analytic family of type (A) with a common domain ${\cal D}$
\[ 
	{\cal D}=\left\{ u\in L^{2}({\bf R}^{3})^{16}\ |\ \nabla_{{\bf y}_{1}}u\in 
	L^{2}({\bf R}^{3})^{16}\right\} .
\] 
Therefore, we note that each eigenvalue $E_{n}(\kappa)$ of $H(\kappa)$
is an analytic function of $\kappa>0$.
Now we claim
\begin{lemma}\label{MS} Let $u\in {\cal H}_{A}^{2}$ be a solution to
the equation (\ref{EV}). Then $u=0$ almost everywhere in ${\bf R}^{6}$.
\end{lemma}
\Proof
We assume that $u$ is nonzero vector in an open subset $\Gamma$ with positive measure.
Then it holds that there exists a nonempty open ball $B_{0}$ in ${\bf R}^{3}$
such that for all ${\bf y}_{2}\in B_{0}$, 
\begin{equation}\label{NOZ} 
	\|Fu(\cdot, {\bf y}_{2})\|_{L^{2}({\bf R}^{3})^{16}}\neq 0.
\end{equation} 
Fix ${\bf y}_{2}\in B_{0}$. Then in view of (\ref{REV}), 
it holds that for some $n$ and every $\kappa$ near $1$
\begin{equation}
	E_{n}(\kappa)=\lambda-k_{0}/|\kappa {\bf y}_{2}|.
\end{equation}
From the analyticity it follows that for any $\kappa>0$
\begin{equation}
	\lambda-k_{0}/|\sqrt{2}\kappa {\bf y}_{2}|=E_{n}(\kappa),
\end{equation} 
which contradicts the fact that $E_{n}(\kappa)\in (-2m,0)\cup (0,2m)$ and
$k_{0}\neq 0$.
\qed

When $m=0$, we see that for any $\kappa>0$, $H(\kappa)$ has no eigenvalues in 
$(-\infty ,0)\cup (0,\infty )$, 
so that we would have
\begin{align}
0=\lambda -\frac{k_{0}}{|\kappa{\bf y}_{2}|}
\end{align}
if eigenvalues existed. It leads to a contradiction.
\qed
\\
{\bf Acknowledgement} The authors would like to thank Professor J. Derezi\'{n}ski,
Professor B. Jeziorski and Doctor M. Oelker for their valuable comments.
The first and the third authors are partially
supported by JSPS Grants-in-Aid No. 22540185 and No.230540226, respectively.
The second and the third authors are indebted to Professors Malcolm Brown, Maria Esteban, 
Karl Michael Schmidt and Heinz Siedentop for the invitation to the programme 
"Spectral Theory of Relativistic Operators" at the Isaac Newton Institute in 
the summer of 2012.
 
\end{document}